%% file: envelope.tex
\documentclass[letterpaper,10pt]{article}
\usepackage[textwidth=11.5cm, textheight=22cm, top=2.0cm, centering]{geometry}
\usepackage{setspace}
\usepackage[T1]{fontenc}
\usepackage{lmodern}
\usepackage[english]{babel}
\usepackage{microtype} 
\usepackage[scaled=0.95]{inconsolata} 
\makeatletter
\renewcommand{\texttt}[1]{{\ttfamily #1}}

\renewcommand{\mathtt}[1]{\text{\texttt{#1}}} 
\makeatother
\usepackage{csquotes}
\usepackage{etoolbox}
\usepackage{bbding}
\usepackage{changepage}
\usepackage{booktabs}
\usepackage{graphicx}
\usepackage[x11names, dvipsnames]{xcolor}
\definecolor{Linkz}{RGB}{30, 110, 170}
\definecolor{Darkenta}{RGB}{185, 35, 90}
\definecolor{Magentz}{RGB}{255, 35, 170}
\definecolor{Lightenta}{RGB}{254, 232, 255}
\definecolor{Reference}{RGB}{35, 180, 90}
\definecolor{Periwinkle}{RGB}{102, 51, 255}
\definecolor{yello}{RGB}{255, 245, 230}
\definecolor{Greeno}{RGB}{0, 140, 100}
\definecolor{Leeno}{RGB}{239, 255, 232}
\definecolor{Nicegreen}{RGB}{100, 200, 130}
\usepackage{amsmath, amsthm, amscd, amssymb, stmaryrd}
\usepackage{stackengine}
\usepackage[square]{natbib}

\usepackage[
colorlinks=true,
linkcolor=Darkenta,
citecolor=Greeno,
urlcolor=Darkenta,
]{hyperref}
\usepackage{titlesec} 
\usepackage{abstract} 
\usepackage{enumitem} 
\usepackage{ccicons} 
\usepackage{float}
\usepackage{placeins} 
\usepackage{tikz} 
\usepackage{tikz-cd}
\usetikzlibrary{positioning,fit,patterns}
\usetikzlibrary{arrows.meta} 
\usepackage{adjustbox} 
\usepackage{pgfmath} 
\usepackage{ifthen} 

\newtheoremstyle{upright}
{6pt plus 2pt minus 2pt} 
{6pt plus 2pt minus 2pt} 
{\normalfont} 
{} 
{\bfseries} 
{.} 
{.5em} 
{} 
\theoremstyle{upright}
\theoremstyle{upright}
\newtheorem{theorem}{Theorem}[subsection]
\newtheorem{remark}[theorem]{Remark}

\newtheorem{definition}[theorem]{Definition}
\newtheorem{proposition}[theorem]{Proposition}
\newtheorem{lemma}[theorem]{Lemma}

\newtheorem{corollary}[theorem]{Corollary}
\newtheorem{example}[theorem]{Example}
\newtheorem{mechanization}[theorem]{Mechanization}
\newtheorem{algorithm}[theorem]{Algorithm}

\makeatletter
\renewenvironment{proof}[1][Proof]{%
	\par\pushQED{\qed}%
	\normalfont
	\topsep6\p@\@plus6\p@\relax
	\trivlist
	\item[\hskip\labelsep\slshape #1\@addpunct{.}]%
}{%
	\popQED\endtrivlist\@endpefalse
}

\makeatother

\usepackage{algpseudocode}
\usepackage[most]{tcolorbox}
\usepackage{listings}
\newtcolorbox{breakbox}[2][]{%
	breakable,
	={#2},
	fonttitle=\bfseries,
	colback=white,
	colframe=black!20,
	coltitle=black,
	colbacktitle=white,
	boxrule=0.5pt,
	arc=0pt,
	boxsep=7pt,
	left=3pt,
	right=2pt,
	top=2pt,
	bottom=4pt,
	fontupper=\small\sffamily, 
	#1
}

\AtBeginEnvironment{quotation}{\sffamily}

\newcommand{\customsectionstyle}[2]{%
	\titleformat{\section}[block]
	{\normalfont\fontsize{#1}{1.2\dimexpr#1\relax}\selectfont\centering}
	{\thesection}{1em}%
	{%
		\ifthenelse{\equal{#2}{true}}{\MakeUppercase}{\relax}%
	}%
}
\newcommand{\customsectionspacing}[3]{%
	\titlespacing*{\section}{#1}{#2}{#3}%
}
\newcommand{\customsubsectionstyle}[2]{%
	\titleformat{\subsection}[block]
	{\normalfont\fontsize{#1}{1.2\dimexpr#1\relax}\selectfont\centering}
	{\thesubsection}{1em}%
	{%
		\ifthenelse{\equal{#2}{true}}{\MakeUppercase}{\relax}%
	}%
}
\newcommand{\customsubsectionspacing}[3]{%
	\titlespacing*{\subsection}{#1}{#2}{#3}%
}
\customsectionstyle{14pt}{true}
\customsectionspacing{0pt}{3.5ex plus 1ex minus 0.2ex}{2.5ex}
\customsubsectionstyle{11pt}{true}
\customsubsectionspacing{0pt}{4ex plus 1ex minus 0.2ex}{2ex}
\usepackage{caption}

\makeatletter
\newcommand{\shorttitle}[1]{\def\@shorttitle{#1}}
\newcommand{\email}[1]{\def\@email{#1}}
\newcommand{\metadata}[1]{\def\@metadata{#1}}
\renewcommand{\maketitle}{%
	\begin{center}
		\vfill
		{\fontsize{18pt}{19pt}\selectfont \@title \par}
		\vspace{1em}
		{\normalsize \@author \par}
		\vspace{0.1em}
		{\normalsize \@date \par}
	\end{center}
}
\makeatother

\begin{document}
\input{content}

\end{document}

%% file: content.tex

\title{\uppercase{Carryless Pairing:\\Additive Pairing\\in the Fibonacci Basis}}
\author{Milan Rosko}
\date{May 2026}

\newcommand{\CL}{\mathsf{CL}}
\newcommand{\Z}{\mathsf{Z}}
\newcommand{\Bdel}{\mathsf{B}}
\newcommand{\pair}{\pi_{\CL}}
\newcommand{\unpair}{\mathsf{unpair}}

\maketitle

\begin{abstract}
	\vspace{0em}
	\footnotesize{
		We define a pairing map $\pair : \mathbb{N}^2\to\mathbb{N}$ that encodes $x$ and $y$ into two disjoint bands of \textsc{Zeckendorf} indices separated by a delimiter computed from $x$. The construction is “carryless” by design: the combined support has no consecutive indices, so each produced code is already in \textsc{Zeckendorf}-normal form, and both evaluation and inversion proceed by additive support operations alone, without multiplication, factorization, or positional digit interleaving. The map is injective not surjective, image membership is decidable by the same support machinery used for decoding. The core correctness theorems are mechanized in \textsc{Rocq}.
	}
\end{abstract}

\section{Introduction}

\subsection{Motivation and origin}
\label{subsec:motivation}

A pair of natural numbers can be encoded using only addition. Working in the Fibonacci (or \textsc{Zeckendorf}) basis, we place the support of $x$ on even Fibonacci indices and the support of $y$ on odd indices beyond a delimiter, and read the resulting set as a single integer.

\begin{figure}[H]
	\centering
	$\pair$
	\vspace{-2ex}
	\includegraphics[width=1\textwidth]{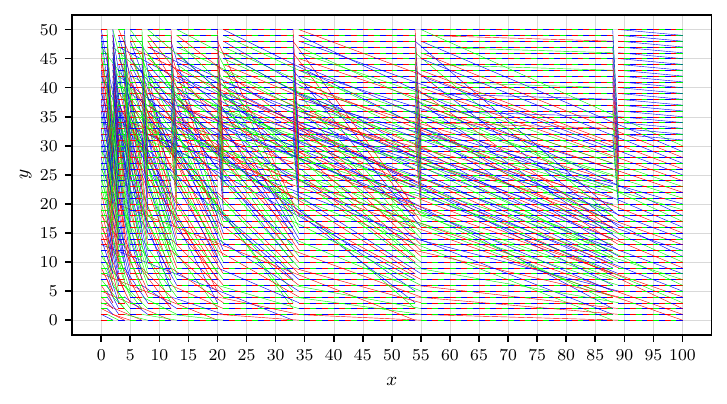}
	\vspace{-2ex}
	\caption{\textsc{Carryless Pairing} provides a fully additive, reversible encoding with a local inverse. The map is injective but not surjective; re-pairing after decoding is a one-step canonicalization onto the image, and image membership is the corresponding fixed-point test.}
	\label{fig:carryless}
\end{figure}

	Our construction belongs to the arithmetic family of pairings rather than to the positional family. The latter schemes---most prominently the \textsc{Z-order} \citep{morton66}---assign digits to fixed positions of a radix expansion. The present construction instead operates on finite \textsc{Zeckendorf} supports and places them into two disjoint bands of Fibonacci indices, separated by a delimiter computed from $x$. Since the bands are disjoint and of opposite parity, the output is already in normal form, and decoding reduces to parity filtering followed by a delimiter test.

	The device emerged from an experiment, namely whether \textsc{Incompleteness} \citep{goedel31,rosser36} can be supported by a pairing scheme whose evaluation step is purely additive in the \textsc{Zeckendorf} basis, dispensing with multiplication and prime factorization, in the spirit of “reverse mathematics” \citep{friedman76}. Its correctness is established under the standard \textsc{Zeckendorf} support interface and validated by a \textsc{Rocq} mechanization (Section~\ref{sec:rocq}). We do not derive that support interface itself from addition-only principles, and conjecture that no such derivation exists.

	The broader motivation is the behaviour of proof obligations under a deliberately weak arithmetic substrate \citep{buss97}. Pairing functions are foundational for provability arguments, yet their image is rarely surveyed in concrete detail. The present construction makes that image, and the band geometry that produces it, directly inspectable.

\section{Background}
\label{sec:background}

\subsection{Pairing and the carryless property}

\begin{definition}[Pairing function]
	In this paper a \emph{pairing function} is a total injective map
	\begin{equation}
	\mathbb{N}^2\to\mathbb{N}
	\end{equation}
	equipped with an effective decoding procedure. Surjectivity is not required; instead we ask for an explicit inversion algorithm and a decidable image test.
\end{definition}

\begin{remark}
	Every natural number admits a non-consecutive \textsc{Zeckendorf} representation, so the available index positions are not exhausted by a single number; there is room for a second number to occupy the unused positions. The construction places all information about $x$ into \emph{even} positions, and all information about $y$ into \emph{odd} positions lying beyond a delimiter $\Bdel(x)$ chosen to the right of every position used by $x$. The two sets of positions are then disjoint and separated, the combined support contains no consecutive indices, and the resulting sum is therefore already in \textsc{Zeckendorf}-normal form---no carry or normalization step is required.
\end{remark}

\begin{definition}[Carrylessness]
	Fix an encoding scheme in which each $n\in\mathbb{N}$ is represented by a finite \emph{support} $S\subseteq I$ and a family of \emph{atoms} $(a_i)_{i\in I}$ with
	\begin{equation}
		n=\sum_{i\in S} a_i,
	\end{equation}
	together with a \emph{normal-form condition} given by a set of forbidden local patterns
	\begin{equation}
	\mathsf{Bad}\subseteq \mathcal{P}_{\mathrm{fin}}(I)
	\quad
	\text{(typically patterns of size $2$ or $3$).}
	\end{equation}
	The support $S$ is \emph{carryless} if it contains no forbidden pattern,
	\begin{equation}
	\forall P\in \mathsf{Bad}\ (P\nsubseteq S);
	\end{equation}
	equivalently, the representation $S$ is already in normal form, and no local rewrite (``carry'' or normalization step) applies. A pairing map $f:\mathbb{N}^2\to\mathbb{N}$ is \emph{carryless} with respect to this scheme if for every $(x,y)$ the code $f(x,y)$ admits a normal-form representation whose support is carryless, so that decoding reads off the components without any prior normalization.
\end{definition}

\begin{definition}[Carryless Fibonacci support]
	The instance relevant to this paper specializes the preceding definition to atoms $a_i = F_i$ and the \textsc{Zeckendorf} normal-form condition $\mathsf{Bad} = \{\{i, i+1\} : i\in\mathbb{N}\}$. A finite support
	\begin{equation}
		\begin{array}{c}
			n = \sum_{e\in S} F_e = F_{e_1}+F_{e_2}+\cdots+F_{e_m},
			\\[0.5em]
			S = \{e_1<e_2<\cdots<e_m\}
		\end{array}
	\end{equation}
	is then “carryless” precisely when it contains no two consecutive indices, in which case the sum is already in \textsc{Zeckendorf}-normal form. A map is “carryless” in this sense if its output support is such on every input.
\end{definition}

\subsection{Contributions and roadmap}

The contributions of the paper are the following:
\begin{enumerate}[label=\textnormal{(\roman*)}]
	\item a carryless, fully additive pairing device over \textsc{Zeckendorf} supports;
	\item explicit pairing and unpairing algorithms relying only on support extraction and bounded scans;
	\item proofs of correctness ($\unpair\circ\pair=\mathrm{id}$) and injectivity;
	\item a one-step image canonicalization test;
	\item a support-level stability lemma and an iterated-coding sketch.
\end{enumerate}

The paper is organized as follows. Section~\ref{sec:background} fixes background notation; Section~\ref{sec:device} defines the device; Section~\ref{sec:algorithms} gives the explicit algorithms; Section~\ref{sec:examples} works through illustrative examples; Section~\ref{sec:correctness} establishes correctness and injectivity; Section~\ref{sec:image} treats image membership, complexity, and stability under local edits; Section~\ref{sec:comparative} contrasts the construction with classical encodings; Section~\ref{sec:rocq} reports the \textsc{Rocq} mechanization.

\subsection{Fibonacci indexing}

We use the standard Fibonacci indexing, cf. \citet{A000045},
\begin{equation}
	F_0 = 0,\quad F_1 = 1,\quad F_{k+2}=F_{k+1}+F_k.
\end{equation}
\textsc{Zeckendorf} supports use indices starting at $2$, so as to avoid the duplication $F_1 = F_2 = 1$. The rank $r(n)$ and intermediate calculations range over the full Fibonacci sequence, including indices $0$ and $1$.

\subsection{\textsc{Zeckendorf} support and rank}

\begin{definition}[Support and Rank]
	For each $n\in\mathbb{N}$, let $\Z(n)\subset\mathbb{N}_{\ge 2}$ be the finite set of indices, all $\ge 2$ and pairwise non-consecutive, such that
	\begin{equation}
		n=\sum_{e\in\Z(n)} F_e.
	\end{equation}
	Define the \emph{rank} as
	\begin{equation}
		r(n)=\min\{e: F_e>n\}.
	\end{equation}
\end{definition}

Existence and uniqueness of $\Z(n)$ are the content of the standard \textsc{Zeckendorf} theorem \citep{zeckendorf,lekkerkerker}. A “greedy” extraction procedure computes $\Z(n)$ by descending Fibonacci indices, at cost $O(r(n))$ when scanning index by index, or $O(|\Z(n)|)$ when admissible indices can be skipped directly. The correctness arguments of Section~\ref{sec:correctness} rely on two standard properties of this support interface: soundness of support summation, and the band-splitting law for carryless encodings.

\section{The carryless pairing device}
\label{sec:device}

\subsection{Carriers and active bands}

\begin{definition}[Carriers and active bands]
	Let $\tau_x=\{e: 2\le e<r(x)\}$ denote the admissible \textsc{Zeckendorf} index range for $x$, and define the \emph{delimiter}
	\begin{equation}
		\Bdel(x)=2\,r(x).
	\end{equation}
	The rank-bounded \emph{even} carrier of $x$ and the rank-bounded \emph{odd} carrier of $y$ relative to $x$ are
	\begin{equation}
		\widehat{\varepsilon}_x=\{2e: 2\le e<r(x)\},\quad
		\widehat{\omega}_{x,y}=\{\Bdel(x)+(2j-1): 2\le j<r(y)\}.
	\end{equation}
	These are bounded sets of admissible indices, fixed from $r(x)$ and $r(y)$ before evaluation; the recursion of the device is confined to them. The corresponding \emph{active bands} are the subsets actually populated by the \textsc{Zeckendorf} supports of $x$ and $y$,
	\begin{equation}
		\begin{array}{c}
				\varepsilon_x^{\mathrm{act}}=\{2e: e\in\Z(x)\}\subseteq \widehat{\varepsilon}_x,\\[1.5ex]
				\omega_{x,y}^{\mathrm{act}}=\{\Bdel(x)+(2j-1): j\in\Z(y)\}\subseteq \widehat{\omega}_{x,y}.
		\end{array}
	\end{equation}
	The odd active band is therefore the image of $\Z(y)$ under $j\mapsto \Bdel(x)+(2j-1)$, not the entire carrier $\widehat{\omega}_{x,y}$.
\end{definition}

\subsection{Definition of carryless pairing}

\begin{definition}[Carryless Pairing]
	Define
	\begin{equation}
		\pair(x,y)=\sum_{k\in \varepsilon_x^{\mathrm{act}}\cup\omega_{x,y}^{\mathrm{act}}} F_k
		=\sum_{e\in\Z(x)} F_{2e}+\sum_{j\in\Z(y)} F_{\Bdel(x)+(2j-1)}.
	\end{equation}
\end{definition}

\begin{lemma}[Carryless Invariant]
	The support of $\pair(x,y)$ contains no adjacent indices.
\end{lemma}

\begin{proof}[Proof sketch.]
	Within the even active band, consecutive elements $e>e'$ of $\Z(x)$ satisfy $e\ge e'+2$ by \textsc{Zeckendorf} validity, hence $2e\ge 2e'+4$. Within the odd active band, consecutive elements $j>j'$ of $\Z(y)$ similarly give
	\begin{equation}
		\Bdel(x)+(2j-1)\ge \Bdel(x)+(2j'-1)+4.
	\end{equation}
	Each band therefore has internal gap at least $4$. Across the delimiter, every populated even index lies below the even number $\Bdel(x)$, while every populated odd index is odd and lies above $\Bdel(x)$; the two bands cannot contain adjacent indices.
\end{proof}

\section{Algorithms}
\label{sec:algorithms}

\begin{algorithm}[Pseudocodes]
	Pairing and unpairing operate on \textsc{Zeckendorf} supports only; the pseudocode below uses the standard support extractor and the parity/delimiter filter from Section~\ref{sec:device}.

	\begin{breakbox}{$\pair$}
	\begin{algorithmic}[1]
		\Procedure{$\pair$}{$x,y$}
			\State $Zx \gets \Z(x)$
			\State $Zy \gets \Z(y)$
			\State $Bx \gets 2\,r(x)$
			\State $S \gets \{2e: e\in Zx\}\cup\{Bx+(2j-1): j\in Zy\}$
			\State \Return $\sum_{k\in S} F_k$
		\EndProcedure
	\end{algorithmic}
	\end{breakbox}

	\begin{breakbox}{$\unpair$}
	\begin{algorithmic}[1]
		\Procedure{$\unpair$}{$n$}
			\State $Zn \gets \Z(n)$
			\State $Zx \gets \{k/2: k\in Zn,\ k\ \text{even}\}$
			\State $x \gets \sum_{e\in Zx} F_e$
			\State $Bx \gets 2\,r(x)$
			\State $Zy \gets \{(k-Bx+1)/2: k\in Zn,\ k\ \text{odd and } k\ge Bx+1\}$
			\State $y \gets \sum_{j\in Zy} F_j$
			\State \Return $(x,y)$
		\EndProcedure
	\end{algorithmic}
	\end{breakbox}

	\begin{breakbox}{$\unpair$ with Check}
	\begin{algorithmic}[1]
		\Procedure{$\unpair_{\mathrm{checked}}$}{$n$}
			\State $(x,y) \gets \Call{$\unpair$}{n}$
			\If{$\pair(x,y)=n$}
				\State \Return $(x,y)$
			\Else
				\State \Return $\bot$
			\EndIf
		\EndProcedure
	\end{algorithmic}
	\end{breakbox}
\end{algorithm}

\section{Examples}
\label{sec:examples}

\begin{example}[Minimal verification]
	At $(x,y)=(1,1)$ we have $\Z(1)=\{2\}$ and $r(1)=3$, hence $\Bdel(1)=6$. The active bands are
	\begin{equation}
		\varepsilon_1^{\mathrm{act}}=\{4\},\quad
		\omega_{1,1}^{\mathrm{act}}=\{9\},
	\end{equation}
	whose combined support is $\{9,4\}$. Summing the corresponding Fibonacci values gives $\pair(1,1)=F_9+F_4=34+3=37$, and indeed $\Z(37)=\{9,4\}$.
\end{example}

\begin{example}[Worked example]
	\label{ex:worked}
	At $(x,y)=(42,1337)$, with $r(42)=10$ and $\Bdel(42)=20$:
	\begin{center}
		\small
		\renewcommand{\arraystretch}{1.8}
		\begin{tabular}{@{}p{1.5cm}p{4.2cm}p{4.2cm}@{}}
			Input & Support & Active band \\
			\hline
			$x=42$ & $\{9,6\}$ & $\{18,12\}$ \\
			$y=1337$ & $\{16,13,11,8,5,3\}$ & $\{51,45,41,35,29,25\}$ \\
			\hline
		\end{tabular}
		\medskip
	\end{center}
	The combined support $\{51,45,41,35,29,25,18,12\}$ has gap at least $4$ throughout, so
	\begin{equation}
		\pair(42,1337) \;=\; 21\,675\,313\,832,
	\end{equation}
	and $\unpair$ recovers $(42,1337)$ by splitting this support at $\Bdel(42)=20$ by parity.
\end{example}

\begin{example}[Outside the image]
	Take $n=6$, with $\Z(6)=\{5,2\}$. Halving the even index gives $x=F_1=1$, hence $\Bdel(1)=6$, and no odd index of $\Z(6)$ lies above $6$, so $y=0$. Since $\pair(1,0)\ne 6$, we conclude $6\notin\mathrm{im}(\pair)$.
\end{example}

\section{Correctness}
\label{sec:correctness}

\begin{theorem}[Inversion]
	For all $x,y\in\mathbb{N}$,
	\begin{equation}
		\unpair(\pair(x,y))=(x,y).
	\end{equation}
\end{theorem}

\begin{proof}[Proof sketch.]
	By the Carryless Invariant, the union $\varepsilon_x^{\mathrm{act}}\cup\omega_{x,y}^{\mathrm{act}}$ is a valid \textsc{Zeckendorf} support, hence equals $\Z(\pair(x,y))$ by uniqueness. Filtering this support for even indices recovers $\{2e:e\in\Z(x)\}$, and halving recovers $\Z(x)$ and thus $x$; the value of $x$ then determines $\Bdel(x)$. Filtering for odd indices above $\Bdel(x)$ recovers $\{\Bdel(x)+(2j-1):j\in\Z(y)\}$, and inverting the affine map recovers $\Z(y)$ and thus $y$.
\end{proof}

\begin{remark}
	The theorem is one-directional. $\unpair$ is a left inverse of $\pair$, but not a right inverse, since $\pair$ is not surjective; the converse identity $\pair(\unpair(n))=n$ characterizes membership in $\mathrm{im}(\pair)$ and is treated in Section~\ref{sec:image}.
\end{remark}

\begin{corollary}[Injectivity]
	$\pair(x,y)=\pair(x',y')$ implies
	\begin{equation}
	(x,y)=(x',y').
	\end{equation}
\end{corollary}

\section{Image}
\label{sec:image}

\subsection{Image membership}

\begin{proposition}[Decidability and one-step canonicalization]
	$n\in\mathrm{im}(\pair)$ if and only if $\pair(\unpair(n))=n$. The test uses only \textsc{Zeckendorf} extraction, parity filtering, and bounded sums, and can be implemented by a bounded greedy scan over $O(r(n))$ indices; whether it is $\Delta_0$ in the sense of bounded arithmetic depends on the chosen formalization of $\Z$, and no formal complexity theorem is certified here.

	Define $\mathsf{can}(n) := \pair(\unpair(n))$. Then $\mathsf{can}$ is idempotent, $\mathsf{can}(\mathsf{can}(n))=\mathsf{can}(n)$; every $n$ reaches an image representative after one application of $\mathsf{can}$, and $n\in\mathrm{im}(\pair)$ iff $\mathsf{can}(n)=n$.
\end{proposition}

\begin{proof}[Proof sketch.]
	If $n\in\mathrm{im}(\pair)$, write $n=\pair(x,y)$; the inversion theorem gives $\unpair(n)=(x,y)$, so $\pair(\unpair(n))=\pair(x,y)=n$. Conversely, $\pair(\unpair(n))=n$ exhibits $n$ as a value of $\pair$. The same argument with $\unpair(n)=(x,y)$ yields $\mathsf{can}(\mathsf{can}(n))=\pair(x,y)=\mathsf{can}(n)$.
\end{proof}

\begin{example}[Canonicalization orbit]
	The standard extractor gives the orbit
	\begin{equation}
		14 \;\xrightarrow{\unpair}\; (1,1) \;\xrightarrow{\pair}\; 37 \;\xrightarrow{\unpair}\; (1,1) \;\xrightarrow{\pair}\; 37 \;\cdots,
	\end{equation}
	so $14\notin\mathrm{im}(\pair)$ and $37=\mathsf{can}(14)$ is its image representative.
\end{example}

\begin{remark}
	No simple closed-form characterization of $\mathrm{im}(\pair)$ is claimed; the natural decision procedure is the support-based fixed-point test above.
\end{remark}

\subsection{Stability and compositionality}
\label{subsec:stability}

\begin{remark}
	After support extraction, pairing costs $O(|\Z(x)|+|\Z(y)|)$ and unpairing $O(|\Z(n)|)$; overall cost is dominated by the support-extraction step, typically $O(r(n))$ for the “greedy” scan. These bounds are informal observations rather than certified theorems.

	By contrast, prime-power encodings (cf.\ \citealp{goedel31}) are $\Delta_0$ only relative to exponentiation, since $n\mapsto p^n$ is not $\Delta_0$-definable in $I\Delta_0$. The carryless device requires no exponentiation: evaluation and inversion use only support extraction and bounded sums.
\end{remark}

\begin{lemma}[Stability under second-coordinate edits]
	Fix $x$. If $\Z(y')=\Z(y)\,\triangle\,\{j\}$, then
	\begin{equation}
		\varepsilon_x^{\mathrm{act}}\cup\omega_{x,y'}^{\mathrm{act}}
		\;=\;
		(\varepsilon_x^{\mathrm{act}}\cup\omega_{x,y}^{\mathrm{act}})\,\triangle\,\{\Bdel(x)+(2j-1)\},
	\end{equation}
	so that
	\begin{equation}
		\pair(x,y')-\pair(x,y)=\pm F_{\Bdel(x)+(2j-1)},
	\end{equation}
	with the sign determined by whether $j$ is added to or removed from $\Z(y)$.
\end{lemma}

\begin{proof}[Proof sketch.]
	The even active band depends only on $x$, hence is unchanged. The odd active band is the image of $\Z(y)$ under the injective map $j\mapsto \Bdel(x)+(2j-1)$, so toggling one source index toggles exactly its image.
\end{proof}

\begin{remark}
	Updates in $x$ are not local in the same sense, since they may shift the delimiter: if $x$ is replaced by $x'$, the odd active band must be recomputed from $\Z(y)$ at the new offset $\Bdel(x')$, and the even active band from $\Z(x')$.
\end{remark}

\begin{remark}
	The binary code iterates by right association. Writing $\pi_{\CL}^{[k]}$ for the $k$-ary iterate,
	\begin{equation}
		\pi_{\CL}^{[k]}(x_1,\ldots,x_k) = \pair(x_1,\pi_{\CL}^{[k-1]}(x_2,\ldots,x_k)),\quad \pi_{\CL}^{[2]} = \pair,
	\end{equation}
	the corresponding round-trip law for the matching right-associated recursive decoder follows by induction from $\unpair\circ\pair=\mathrm{id}$. The present paper does not develop a separate $k$-ary interface.
\end{remark}

\section{Applications}

We highlight four directions in which the device is naturally applied.

\begin{enumerate}[label=\textnormal{(\roman*)}]
	\item \emph{Additive consistency arguments.} The construction provides a concrete additive mechanism whose correctness yields consistency-style facts usable in arguments that avoid multiplication (cf.\ \citealp{buss97,pudlak2008}).
	\item \emph{Automated proof checking.} The separation of supports lends itself to a small, checkable Hilbert system with fast verification (see \citealp{kleene1952,pudlak98} for proof-theoretic context).
	\item \emph{Image structure.} The image of $\pair$ exhibits substantial irregularity and is itself an object of interest (see Figure~\ref{fig:carryless}).
	\item \emph{Reverse mathematics and coding strength.} In weak arithmetics, uniform coding of tuples is not a mere convenience but can mark a genuine increase in expressive power (cf.\ \citealp{friedman76,simpson09}); a natural question is which comprehension or induction axioms are required to derive $\unpair\circ\pair=\mathrm{id}$ and its iterated forms.
\end{enumerate}
\section{Comparative encodings}
\label{sec:comparative}

We distinguish two families of encodings: \emph{arithmetic} encodings, in which structure is carried by number-theoretic operations, and \emph{positional} encodings, in which structure is carried by digit positions. The standard arithmetic baselines are the multiplicative (\textsc{Gödel}) and polynomial (\textsc{Cantor}) pairings; the standard positional baseline is \textsc{Bit-interleaving}.

\begin{figure}[H]
	\centering
	$\pi_{\mathrm{GN}}$
	\vspace{-2ex}
	\includegraphics[width=1\textwidth]{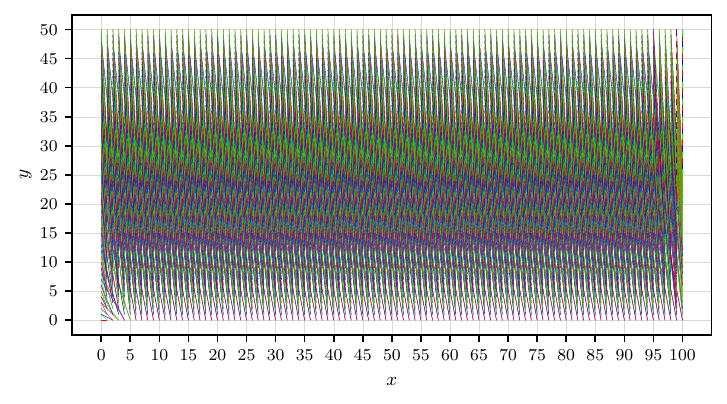}
	\vspace{-2ex}
	\caption{\textsc{Gödel numbering} encodes sequences via prime powers; inversion proceeds by divisibility and factorization.}
	\label{fig:godel}
\end{figure}

The classical \textsc{Gödel numbering} represents a finite sequence by a product of prime powers and decodes structure through divisibility and factorization (Fig.~\ref{fig:godel}). It presupposes multiplication and exponentiation; the cost contrast with the carryless device was made precise in Section~\ref{sec:image}.

\begin{figure}[H]
	\centering
	$\pi_{\mathrm{C}}$
	\vspace{-2ex}
	\includegraphics[width=1\textwidth]{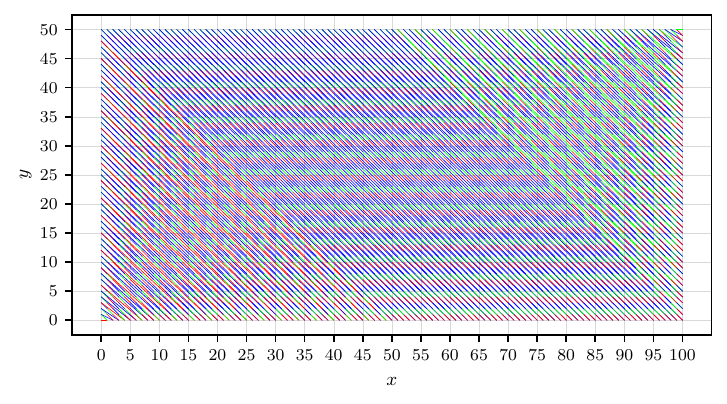}
	\vspace{-2ex}
	\caption{\textsc{Cantor pairing} is a quadratic bijection on $\mathbb{N}^2$.}
	\label{fig:cantor}
\end{figure}
The \textsc{Cantor pairing} function (Fig.~\ref{fig:cantor}) is the polynomial bijection
\begin{equation}
	\pi_{\mathrm{C}}(x,y)=\tfrac{1}{2}(x+y)(x+y+1)+y,
\end{equation}
whose inversion is closed-form but uses multiplication and an integer square root. Unlike the carryless device, it is bijective; unlike the carryless device, it requires arithmetic beyond addition.

\begin{figure}[H]
	\centering
	$\pi_{\mathrm{B}}$
	\vspace{-2ex}
	\includegraphics[width=1\textwidth]{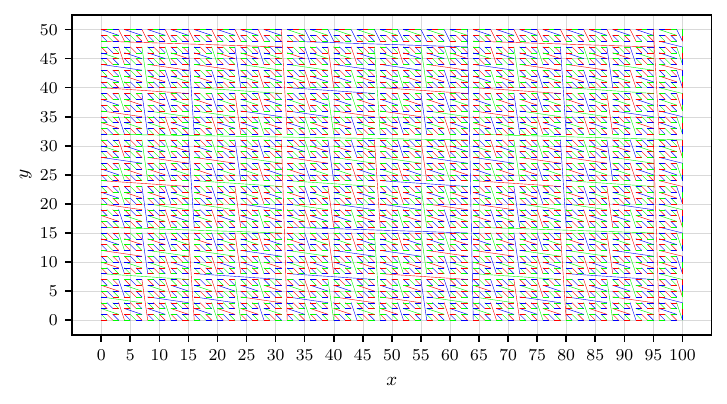}
	\vspace{-2ex}
	\caption{\textsc{Bit-interleaving} yields the \textsc{Z-order} used for multidimensional indexing.}
	\label{fig:morton}
\end{figure}
The \textsc{Z-order} of \citet{morton66} interleaves the binary digits of $(x,y)$ to obtain a total, surjective map (Fig.~\ref{fig:morton}). The bit positions are fixed by parity,
\begin{equation}
	\pi^{\mathrm{B}}(x,y):\quad \mathrm{pos}(b^x_e)=2e,\quad \mathrm{pos}(b^y_j)=2j+1,
\end{equation}
so the placement of each digit is independent of the other input. The carryless device differs in a complementary way: the $x$-band sits at fixed even positions $2e$, but the $y$-band is shifted past the rank-determined delimiter $\Bdel(x)=2\,r(x)$. Updating $x$ does not redistribute the $y$-support digit by digit; it applies a single affine translation $j\mapsto\Bdel(x)+(2j-1)$ to the image of $\Z(y)$.

\section{Mechanization}
\label{sec:rocq}

\begin{mechanization}
The construction and its correctness are formalized in \textsc{Gallina}, the specification language of:

\begin{center}
	\vspace{0.1em}
	\includegraphics[width=3cm]{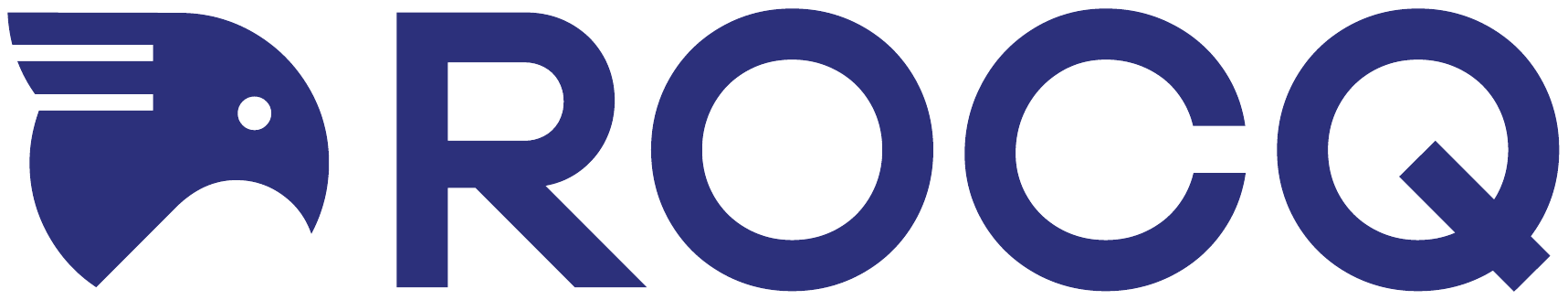}\\
	\vspace{0.1em}
\end{center}

	The paper's $\pair$ and $\unpair$ appear in the source as \texttt{encode} and \texttt{decode}. The two headline theorems are \texttt{decode\_encode}, which establishes $\unpair\circ\pair=\mathrm{id}$, and \texttt{encode\_injective}, which establishes injectivity. Both are closed under the global context---no axioms, no admits.

\begin{center}
	\href{https://github.com/Milan-Rosko/Proofcase/tree/main/theories/A001}{Github.com $\to$ Milan-Rosko $\to$ Proofcase $\to$ A001}
\end{center}

	The development is layered to mirror the logical structure of the paper. The Fibonacci basis and the “greedy” extractor are introduced first, and the support interface assumed informally in Section~\ref{sec:background} is itself proved: existence and uniqueness of the support are the certified theorems \texttt{Z0\_sound}, \texttt{Z0\_valid}, and \texttt{Zeckendorf\_unique}, so the support primitives the paper relies on are not axioms of the \textsc{Rocq} development.

	The carryless decomposition of Section~\ref{sec:device} is then established as three named theorems
	\begin{equation}
		\begin{array}{l}
			\mathtt{Z0\_pair\_is\_concat},\\
			\mathtt{Z0\_even\_split},\\
			\mathtt{Z0\_odd\_split},
		\end{array}
	\end{equation}
	which witness the parity split used in the proof of inversion.

	The repository exposes an executable surface for the construction, available both in \textsc{Gallina} and as extracted \textsc{OCaml}. The diagnostic
	\begin{equation}
		\mathtt{Compute\_Pair\_Unpair\_Check}: \mathbb{N} \times \mathbb{N} \to \mathbb{N} \times (\mathbb{N} \times \mathbb{N}) \times \mathbb{N} \times \mathtt{Code\_Status}
	\end{equation}
	runs the full $\pair$/$\unpair$/canonicalization loop on $(a,b)$ in one call, returning the code, its decoding, its re-encoding, and the canonicality status. Its components are individually exposed:
	\begin{equation}
		\begin{array}{l}
			\mathtt{Check\_Pairing}: \mathbb{N} \to \mathbb{N},\\
			\mathtt{In\_Imageb}: \mathbb{N} \to \mathtt{bool},\\
			\mathtt{Status\_Of\_Code}: \mathbb{N} \to \mathtt{Code\_Status},
		\end{array}
	\end{equation}
	computing respectively the canonicalization $\mathsf{can}(n)$, the fixed-point test deciding $n\in\mathrm{im}(\pair)$, and a tagged form of that test. The structured dispatcher
	\begin{equation}
		\mathtt{A001\_IO}: \mathtt{IO\_Query} \to \mathtt{IO\_Result}
	\end{equation}
	packages pairing and inspection queries under a single entry point.
\end{mechanization}

\clearpage

\vspace*{\fill}

\section*{References and Notes}

{\scriptsize
	\bibliographystyle{plainnat}
	\setlength{\bibsep}{0.1em}
	\bibliography{refs}}

\subsection*{Version History}

	The first manuscript of this line, from September 2025, framed the \textsc{Carryless Pairing} as a tool inside an additive reformulation of \textsc{Incompleteness} with applications to diophantine questions. As the construction stabilized, it became apparent that the pairing function merits independent technical treatment, and the foundational framing was relegated to companion work. The present, fifth manuscript refines the fourth version.

\subsection*{Final Remarks}

	The author welcomes criticism, proposed extensions, scholarly correspondence and constructive dialogue. No conflicts of interest are declared. This research received no funding.
\begin{center}
	\scriptsize{
	\vspace{3em}

	Milan Rosko

	\vspace{1em}
	ORCID: \href{https://orcid.org/0009-0003-1363-7158}{\textsf{0009-0003-1363-7158}}\\
	Email: \href{mailto:Q1012878@studium.fernuni-hagen.de}{\textsf{Q1012878@studium.fernuni-hagen.de}}\\
	Email: \href{mailto:hi-at-milanrosko.com}{\textsf{hi-at-milanrosko.com}}\\
	\vspace{1em}
	Licensed under \enquote{Deed} \ccby \\ \href{http://creativecommons.org/licenses/by/4.0/}{\scriptsize\textsf{creativecommons.org/licenses/by/4.0}}
	}
\end{center}

\vspace*{\fill}

%% file: envelope.bbl
\begin{thebibliography}{12}
\providecommand{\natexlab}[1]{#1}
\providecommand{\url}[1]{\texttt{#1}}
\expandafter\ifx\csname urlstyle\endcsname\relax
  \providecommand{\doi}[1]{doi: #1}\else
  \providecommand{\doi}{doi: \begingroup \urlstyle{rm}\Url}\fi

\bibitem[Buss(1997)]{buss97}
S.~R. Buss.
\newblock {Bounded Arithmetic, Cryptography and Complexity}.
\newblock \emph{Theoria}, 63\penalty0 (3):\penalty0 147--167, 1997.
\newblock URL \url{https://doi.org/10.1111/j.1755-2567.1997.tb00745.x}.

\bibitem[Friedman(1976)]{friedman76}
H.~M. Friedman.
\newblock {Systems of Second-Order Arithmetic with Restricted Induction, I,
  II}.
\newblock \emph{Journal of Symbolic Logic}, 41\penalty0 (2):\penalty0 557--559,
  1976.
\newblock URL \url{https://doi.org/10.2307/2272259}.
\newblock Meeting of the Association for Symbolic Logic.

\bibitem[G{\"o}del(1931)]{goedel31}
K.~G{\"o}del.
\newblock {{\"U}ber formal unentscheidbare S{\"a}tze der Principia Mathematica
  und verwandter Systeme I}.
\newblock \emph{Monatshefte f{\"u}r Mathematik und Physik}, 38\penalty0
  (1):\penalty0 173--198, 1931.
\newblock URL \url{http://doi.org/10.1007/BF01700692}.

\bibitem[Kleene(1952)]{kleene1952}
S.~C. Kleene.
\newblock \emph{{Introduction to Metamathematics}}.
\newblock North-Holland, 1952.
\newblock ISBN 9780444896230.

\bibitem[Lekkerkerker(1952)]{lekkerkerker}
C.~G. Lekkerkerker.
\newblock {Representatie van natuurlijke getallen door een som van getallen van
  {Fibonacci}}.
\newblock \emph{Simon Stevin}, 29:\penalty0 190--195, 1952.

\bibitem[Morton(1966)]{morton66}
G.~M. Morton.
\newblock {A computer Oriented Geodetic Data Base; and a New Technique in File
  Sequencing}.
\newblock 1966.
\newblock IBM Ltd. Technical Report.

\bibitem[Pudlák(1998)]{pudlak98}
P.~Pudlák.
\newblock \emph{{The Lengths of Proofs}}, volume~28 of \emph{Studies in Logic
  and the Foundations of Mathematics}.
\newblock Springer, 1998.
\newblock ISBN 9780444898401.

\bibitem[Pudlák(2008)]{pudlak2008}
P.~Pudlák.
\newblock {Fragments of bounded arithmetic and the lengths of proofs}.
\newblock \emph{Journal of Symbolic Logic}, 73\penalty0 (4):\penalty0
  1389–1406, 2008.
\newblock URL \url{http://doi.org/10.2178/jsl/1230396927}.

\bibitem[Rosser(1936)]{rosser36}
J.~B. Rosser.
\newblock {Extensions of Some Theorems of {G{\"o}del} and {Church}}.
\newblock \emph{Journal of Symbolic Logic}, 1\penalty0 (3):\penalty0 87--91,
  1936.
\newblock URL \url{https://doi.org/10.2307/2269028}.

\bibitem[Simpson(2010)]{simpson09}
S.~G. Simpson.
\newblock \emph{{Subsystems of Second Order Arithmetic}}.
\newblock Perspectives in Logic. Cambridge University Press, 2nd edition, 2010.
\newblock URL \url{https://doi.org/10.1017/CBO9780511581007}.

\bibitem[Sloane(1964)]{A000045}
N.~J.~A. Sloane.
\newblock {A000045: Fibonacci numbers}.
\newblock The On-Line Encyclopedia of Integer Sequences, 1964.
\newblock URL \url{https://oeis.org/A000045}.

\bibitem[Zeckendorf(1972)]{zeckendorf}
E.~Zeckendorf.
\newblock {Repr{\'e}sentation des nombres naturels par une somme de nombres de
  {Fibonacci} ou de nombres de {Lucas}}.
\newblock \emph{Bulletin de la Soci{\'e}t{\'e} Royale des Sciences de
  Li{\`e}ge}, 41\penalty0 (4-6):\penalty0 179--182, 1972.

\end{thebibliography}
